\documentclass[review]{elsarticle}%

\journal{Journal of \LaTeX\ Templates}

\usepackage{amssymb}
\usepackage{colortbl}
\usepackage{amsmath}

\usepackage{moreverb}
\usepackage{booktabs}
\usepackage{graphicx}

\usepackage{amsthm}

\newcommand{\diag}{\mathop{\textrm{diag}}}

\newtheorem{thm}{Theorem}
\newtheorem{lem}{Lemma}
\newtheorem{rem}{Remark}

\newtheorem{defi}{Definition}

\journal{}

\begin{document}

\begin{frontmatter}

\title{Multiple integral inequalities and stability analysis of time delay systems}

\author{\'E. Gyurkovics \footnotemark[1], 
T. Tak\'acs\footnotemark[2]}

\address{\footnotemark[1] Mathematical Institute, Budapest University of Technology and Economics, Budapest, Pf. 91,
1521, Hungary \\
 \footnotemark[2] Corvinus University of Budapest,
8 F\H{o}v\'am t\'er, H-1093, Budapest, Hungary}

\begin{abstract}
This paper is devoted to stability analysis of continuous-time delay systems based on a set of Lyapunov-Krasovskii functionals. New multiple integral inequalities are derived that involve the famous Jensen's and Wirtinger's inequalities, as well as the recently presented
Bessel-Legendre inequalities of A. Seuret and F. Gouaisbaut, (2015) \cite{seur14b},
and the Wirtinger-based multiple-integral inequalities of M. Park et al. (2015) and T.H. Lee et al. (2015)
\cite{park15, leejfi15}.
 The present paper aims at showing that the proposed set of sufficient stability conditions can be arranged into a bidirectional hierarchy of LMIs establishing a rigorous theoretical basis for comparison of conservatism of the investigated methods. Numerical examples illustrate the efficiency of the method.
\end{abstract}

\begin{keyword}
Integral inequalities, stability analysis, continuous-time delay systems, hierarchy of LMIs
\end{keyword}

\end{frontmatter}

\section{Introduction}
Time delays are present in many physical, industrial and engineering
systems. The delays may cause instability or poor performance of systems, therefore much attention
has been devoted to obtain tractable stability criteria of systems with time delay during the past few decades
(see e.g. the monographs \cite{Briat14}-
\cite{WHSh}, some recent papers \cite{seur14b}, \cite{park15}, \cite{kim16}-\cite{ze15} 
 and the references therein).
Several approaches have been elaborated and successfully applied for the stability analysis of time delay systems (see the references above for excellent overviews).

Lyapunov method is one of the most fruitful fields in the stability analysis of time delay systems. On the one hand, more and more involved Lyapunoov-Krasovskii functionals (LKF) have been introduced during the past decades. On the other hand, much effort has been devoted to derive more and more tight inequalities (Jensen's inequlity and different forms of Wirtinger's inequality \cite{seur14b}-
\cite{frid2014}, \cite{kim16}-\cite{sga15},
 \cite{gye15}, \cite{zha15}, etc.) for the estimation of quadratic single, double and multiple integral terms in the derivative of the LKF. Simultaneously, augmented state vectors are introduced in part as a consequence of the improved estimations, in part in an ad'hoc manner. The effectiveness of different methods is mainly compared using some numerical examples. Recently, the authors of \cite{seur14b}, \cite{seur14} have introduced a very appealing idea of the hierarchy of LMI conditions offering a rigorous theoretical basis for comparison of stability LMI conditions. Based on Legendre polynomials, they proposed a generic set of single integral inequalities opening the way to the derivation of a set of stability conditions forming a hierarchy of LMIs. A further possibility for the derivation of improved stability conditions have been proposed by \cite{park15} and \cite{leejfi15} using multiple integral quadratic terms in the LKF, together with Wirtinger-based multiple integral inequalities. Naturally the question arises: how these two lines of investigations are related to each other, and how sufficient stability conditions can be derived unifying the approaches of using multiple integral quadratic terms in the LKF and refined estimations of these integral terms.

\textit{The aim of the present work is to answer these questions. On the one hand,
multiple integral inequalities based on orthogonal hypergeometric polynomials will be derived that extend the results of \cite{seur14b}, \cite{seur14} to multiple integrals and improve the estimations of \cite{park15} and \cite{leejfi15}. On the other hand, a multi-parametric set of LMI conditions will be constructed, and it will be shown that a two parametric subset forms a bidirectional hierarchy of LMIs.
}

Analogous results have been presented for discrete-time systems in \cite{gytkn}.

The paper is organized as follows. In Section 2
it is shown, how the quadratic terms of the derivative of the LKF can be estimated by Bessel-type
inequalities.  It is also proven that these estimations relevantly improve a
recently published result. A sufficient condition of asymptotic stability is presented in the form of an LMI in Section 3.
The hierarchy of LMI conditions is established then in Section 4. Some benchmark numerical examples are shown
in Section 5, the results of which are compared to earlier ones known from the literature.
Finally, the conclusions will be drawn.

The notations applied in the paper are very standard, therefore we mention only a few of them. 
 Symbol $A \otimes B$ denotes the Kronecker-product of matrices $A,B,$ while $\mathbf{S}_{n}$ and $\mathbf{S}_{n}^{+}$ are the set of symmetric and positive definite symmetric matrices of size $n\times n,$ respectively. 

\section{Multiple integral inequalities}

\subsection{Preliminaries}

The paper deals with the stability analysis of the following continuous-time time delay system
\begin{eqnarray}
\dot{x} (t)&=&Ax(t)+A_{d_1} x(t-\tau )
+ A_{d_2}\int_{t-\tau}^{t} x(s) ds, \; t \geq 0, \label{x1} \\
x_0 (t) &=& \varphi (t), \; t \in [-\tau,0],      \label{xx1}
\end{eqnarray}
where $x(t) \in \mathbf{R}^{n_x}$ is the state, $A$, $A_{d_1}$and $A_{d_2}$ are given constant matrices of
appropriate size, the time delay $\tau $ is a known positive integer and $x_0(.)$ is the initial function.

\vspace{0.3cm}
A.) \textit{A Bessel-type inequality. }
Let $\mathbf{E}$ be a Euclidean space with the scalar product $\langle . , . \rangle$, and let $\pi _i \in \mathbf{E}, \; (i=0,1, \ldots )$ form an orthogonal system. For any $f,g \in \mathbf{E}^{n},$ define $ \langle f, g\rangle =\sum _{i=1}^{n} \langle f_i, g_i\rangle . $
%
Let $W \in \mathbf{S}^{+}_{n}$. For any $f\in\mathbf{E}^{n},$ consider the functional
\begin{eqnarray}
J_{W}(f) =  \langle f,W f \rangle .   \label{e4}
\end{eqnarray}
\begin{lem} \label{lem:11}
If $\nu \geq 0$ is a given integer,
then the following inequality holds
\begin{eqnarray}
J_{W} (f)\geq \sum _{j=0}^{\nu}\frac{1}{\left\| \pi_{j} \right\|^2} w_{j}^{T}   W  w_j,     \label{e6}
\end{eqnarray}
where $w_j = \langle f, \pi_{j}\rangle ,$ and the scalar product is taken componentwise.
\end{lem}

\textbf{Proof.} The proof is standard, therefore it is omitted.
%

\vspace{0.3cm}
B.) \textit{Orthogonal hypergeometric polynomials. }
Suppose that $m\geq0$ is a given integer and consider the closed interval $[a,b].$
For functions $g_1,g_2 \in L_2 [a,b] $ define a scalar product by
\begin{eqnarray}
\langle g_1,g_2\rangle _{m,[a,b]} = \int _{a}^{b} \left( \frac{s-a}{b-a} \right)^{m} g_1 (s) g_2 (s) ds.  \label{e2}
\end{eqnarray}
It is easy to see that $\langle g_1,g_2\rangle _{m,[a,b]}$ can equivalently be expressed  as
\begin{eqnarray}
\langle g_1,g_2\rangle _{m,[a,b]} &=& \frac{m!}{(b-a)^m} \int _{a}^{b} \int _{v_1}^{b} ... \int _{v_m}^{b} g_1 (s) g_2 (s)ds dv_m ... dv_1
, \mbox{ if }  m>0. \label{e1}
\end{eqnarray}
 (If $m=0,$ then a single integral is considered.)
Substitute $s \in [a,b]$ by $s=a+(b-a)x$, where $x \in [0,1],$ and set  $G_i (x)=g_i (a+(b-a)x),$ $(i=1,2)$ on the right hand side of (\ref{e2}),
then we obtain that
\begin{eqnarray}
\langle g_1,g_2\rangle _{m,[a,b]} = (b-a) \int _{a}^{b} x^{m} G_1 (x) G_2 (x) dx = (b-a) \langle G_1,G_2\rangle _{m,[0,1]}. \label{gy30}
\end{eqnarray}
Thus it is sufficient to consider the orthogonal polynomials with respect to $\langle .,.\rangle  _{m,[0,1]}.$

For any fixed non-negative integer $m$, let us denote by $P_{m,n}, \; (n=0,1, \ldots)$ the polynomials of degree $n$ orthogonal with respect to $\langle .,.\rangle  _{m,[0,1]}.$  (For general theory see e.g. \cite{gau}.)
They can be given by the two parameters generalization of the Rodrigues-formula:
\begin{eqnarray}
P_{m,0}(x) & \equiv& 1, \label{ee1}\\
P_{m,n}(x) &=& \frac{1}{n!} \frac{1}{x^m} \frac{d^n}{dx^n} \left( x^m (x^2-x)^n \right) , \hspace{0.5cm} n=1,2,... \label{e3}
\end{eqnarray}
For $m=0,$ this is the usual Rodrigues formula for the shifted Legendre polynomials.

We note that that polynomials (\ref{ee1})-(\ref{e3}) satisfy certain hypergeometric-type differential equation (see e.g. \cite{Area03} and \cite{San97}). This is why they are frequently called "orthogonal hypergeometric polynomials". By straightforward calculation, it can be shown that they have the properties
\begin{eqnarray}
&(i) \;\;& \left\| P_{m,n} \right\| _{m,[0,1]}^{2} = \int  _{0}^{1} x^m P_{m,n}^{2} (x) dx = \frac{1}{m+2n+1} ,    \label{e3b} \\
&(ii) \;&P_{m,n}(0)=(-1)^n  \frac{m+n}{n} , \hspace{1cm} P_{m,n}(1)=1.
\end{eqnarray}
The polynomials
\begin{eqnarray}
p_{m,n} (t) = P_{m,n} \left( \frac{t-a}{b-a} \right)   \label{e3a} 
\end{eqnarray}
are orthogonal with the scalar product (\ref{e2}), and
\begin{equation}\label{gy01}
  \left\| p_{m,n} \right\| _{m,[a,b]}^{2}=\frac{b-a}{m+2n+1}, \hspace{0.5cm} p_{m,n}(a)=(-1)^n  \frac{m+n}{n} ,
  \hspace{0.5cm}p_{m,n}(b)=1.
\end{equation}

\subsection{Integral inequalities }

Let $W \in \mathbf{S}^{+}_{n},$ $[a,b] \subset \mathbf{R}$ with $b-a>0$ and $0 \leq m \in \mathbf{Z}$ be given. For any continuous $f: [a,b] \rightarrow \mathbf{R}^{n},$ consider the functional
\begin{eqnarray}
J_{W,m,a,b} (f) = \frac{m!}{(b-a)^m} \int _{a}^{b} \int _{v_1}^{b} ... \int _{v_m}^{b} f^{T}(s)W f(s)ds dv_m ... dv_1,   \label{e400}
\end{eqnarray}
which can also be expressed as    
\begin{eqnarray}
J_{W,m,a,b} (f)  = \int _{a}^{b} \left( \frac{s-a}{b-a}  \right)^{m} f^T (s) W f(s) ds = \langle f, W f \rangle _{m,[a,b]}. \label{e5}
\end{eqnarray}
Let  
$\nu _m\geq0$ be a given  integer. 
One can apply now Lemma \ref{lem:11} with $\mathbf{E}=L_{[a,b]}^2,$ the scalar product (\ref{e2}), $\nu =\nu _m$ and   $\pi _j=p_{m,j}, $ $(j=0,1,\ldots, \nu _m).$
Now, \textit{our aim is to derive a lower estimation as a quadratic form with respect to variables independent of $m.$ }


\begin{lem} \label{lem:21}
 Let $M>0$ and $\nu _m\geq0$ be  given  integers satisfying the condition $m+\nu _m \leq M-1$.
Let $J_{W,m,a,b} (f)$ be defined by (\ref{e5}). Then the following inequality holds true:
\begin{eqnarray}
J_{W,m,a,b} (f) \geq \frac{1}{b-a} \Phi_ M^{T} \left( \Xi_m \otimes I \right)^T \mathcal{W}_m \left( \Xi_m \otimes I \right) \Phi_ M, \label{e10}
\end{eqnarray}
where $\mathcal{W}_m = \mbox{\emph{diag}} \left\{ (m+1), (m+3), \ldots, (m+2\nu_m +1) \right\}\otimes W, $
$\Phi_ M^T =
\begin{bmatrix}\phi_0^T,&\ldots ,&\phi_{M-1}^T\end{bmatrix} \; $ with
$\phi _l = \int _{a}^{b} p_{0,l} (s)f(s)ds,$
and matrix $\Xi_m$ is given by (\ref{gy20}) below.
\end{lem}
\textbf{Proof.}
 Introduce the notation $w_{m,j}=\langle f,p_{m,j}\rangle _{m,[a,b]}$  $(j=0,1,\ldots,\nu _m, )$  needed to apply Lemma \ref{lem:11}.
Clearly,
\begin{eqnarray}
w_{m,j} =  \int _{a}^{b}
\left( \frac{s-a}{b-a}\right)^m P_{m,j} \left( \frac{s-a}{b-a} \right) f(s) ds .  \label{e8}
\end{eqnarray}
The degree of polynomials $q_{m,m+j}(x)=x^m P_{m,j} (x)$ appearing in (\ref{e8})
is exactly $m+j,$ thus these polynomials can be expressed as
\begin{eqnarray}
q_{m,m+j} (x) = \sum _{l=0}^{M-1} \xi _{j,l}^{m} P_{0,l} (x),   \label{e9}
\end{eqnarray}
where $\xi_{j,l}^{m}=0,$ if $m+j<l \leq M-1.$  Using the definition of $\phi _l$ and $\Phi_ M$ we obtain
\begin{eqnarray}
w_{m,j}  = \sum _{l=0}^{M-1} \xi _{j,l}^{m} \phi _l = \left(\underline{\xi}_{j}^{m}\otimes I\right) \Phi _M,  \label{e900}
\end{eqnarray}
where $\underline{\xi}_{j}^{m}=\left( \xi_{j,0}^{m}, \ldots, \xi_{j,M-1}^{m} \right).$
Introduce the notation
\begin{equation}\label{gy20}
  \Xi_m =
\left[\left(\underline{\xi}_{0}^{m}\right)^T, \ldots,\left(\underline{\xi}_{\nu_m}^{m}\right)^T \right]^T \in \mathbf{R}^{(\nu_m +1)\times M}
\end{equation}
Estimation (\ref{e10}) can be obtained by direct substitution taking into account (\ref{gy30}) and (\ref{e3b}).  $\Box$
\begin{rem}\label{rem:1}
We note that Lemma \ref{lem:21} is closely connected with Theorem 2.2 of \cite{zha15}. Both results are based (explicitly or implicitly) on Lemma \ref{lem:11}, thus they are substantially equivalent. The estimation of Lemma \ref{lem:21} may be more advantageous when it is applied for derivative of functions (see Lemma \ref{lem:31} below) and  for stability analysis of time delay systems. The advantage is twofold: on the one hand, the variables are expressed using a common set of orthogonal polynomials independent of $m,$ on the other hand, the dependence on the length of the interval is relatively simple, since the matrices $\Xi _m$ and $\mathcal{W}_m$ do not depend on $b-a.$
\end{rem}
\begin{rem}\label{rem:2}
Paper \cite{zha15} gives a thorough and detailed discussion of the relation between their WOPs-based result and the Jensen's and Wirtinger's inequalities published in a wide range of previous literature, therefore we only compare Lemma \ref{lem:21} to the recently published multiple integral inequality of Lemma 5 of \cite{leejfi15}. Using the notations of \cite{leejfi15}, we can see that the relation of the investigated functionals can be given as
\begin{equation}\label{LP1}
  G_l(f,a,b,W)=\frac{(b-a)^l}{l!}J_{W,m,a,b} (f).
\end{equation}
To express the estimation of the present paper with the variables of \cite{leejfi15}, we need a short computation to show that
$w_{l,0}=\frac{l!}{(b-a)^l}g_l(f,a,b),$ and  $w_{l,1}=-\frac{(l+1)!}{(b-a)^l}\Upsilon_l(f,a,b).$ Employing the proposed estimation with (\ref{LP1}), we obtain that
\begin{eqnarray}
   G_l(f,a,b,W) &=& 
   \geq \frac{(l+1)!}{(b-a)^{l+1} } g_l^T(f,a,b)W g_l(f,a,b) \nonumber\\
  &  & \hspace{0.5cm}+ (l+1)^2\frac{l!(l+3)}{(b-a)^{l+1}}\Upsilon_l^T(f,a,b)W\Upsilon_l(f,a,b).
   \label{LP2}
\end{eqnarray}
The first term of the lower bound of \cite{leejfi15} is the same as in (\ref{LP2}), while the second term is smaller inasmuch as it has the coefficient $1$ in place of $(l+1)^2,$ thus the estimation of the present paper is tighter.

The considerations above indicate that the results of \cite{leejfi15} correspond to the choice of $\nu _l =1,$ but
the authors of \cite{park15} do not derive any estimation that can be characterized with $\nu _l >1.$
\end{rem}


Next, we shall derive a lower estimation also for the case, when the functional is applied to the derivative $f'(s)=\frac{d}{ds}f(s),$ 
i.e. consider 
\begin{eqnarray}
J_{W,m,a,b} (f^{\prime}) &=& 
 \int _{a}^{b} \left( \frac{s-a}{b-a}  \right)^{m} f^{\prime}(s)^T W f^{\prime}(s) ds  
 = \langle f^{\prime}, W f^{\prime} \rangle _{m,[a,b]}. \label{gy40}
\end{eqnarray}
\begin{lem} \label{lem:31}
Let $\ M \ $  and $\ \nu _m \ $ be  given non-negative integers satisfying the condition  $\ m+\nu _m
\leq \; $ $\max \left\{0,M-1\right\}$.
Let $J_{W,m,a,b} (f')$ be defined by (\ref{gy40}). Then the following inequality holds true:
\begin{eqnarray}
J_{W,m,a,b} (f^{\prime}) \geq \frac{1}{b-a} \widetilde{\Phi}_ M ^{T} \left( \mathcal{Z}_m \otimes I \right)^T \mathcal{W}_m
\left( \mathcal{Z}_m \otimes I \right) \widetilde{\Phi}_ M, \label{e12}
\end{eqnarray}
where $\mathcal{W}_m $ is the same as in Lemma \ref{lem:21}, $\widetilde{\Phi}_ M = \mbox{col} \left\{ f(b), \; f(a), \; \frac{1}{b-a}\phi_0, \ldots , \frac{1}{b-a}\phi_{M-1} \right\},$  if $ M>0$,
$\widetilde{\Phi}_ 0 = \mbox{col} \left\{ f(b), \; f(a) \right\},$
 and matrix $\mathcal{Z}_m $ is given by (\ref{gy50}) below.
\end{lem}
\textbf{Proof.}  Set $\theta_{m,j} = \langle f^{\prime},p_{m,j}\rangle _{m,[a,b]},$ $(j=0,1,\ldots,\nu _m, )$ and apply Lemma \ref{lem:11}. In order to
obtain the estimation (\ref{e12}), we have to perform a short computation, as follows.
Consider the previously introduced polynomials $q_{m,m+j}$ again,  then integrating by parts we obtain
\begin{eqnarray}
\theta_{m,j} &=& \int _{a}^{b} \left( \frac{s-a}{b-a}   \right)^m p_{m,j} (s) f^{\prime }(s) ds =  \nonumber   \\
&=& q_{m,m+j} (1)f(b) - q_{m,m+j} (0) f(0) - \frac{1}{b-a} \int _{a}^{b} q^{\prime }_{m,m+j} \left( \frac{s-a}{b-a} \right) f(s)ds.
  \label{e11}
\end{eqnarray}
Express now polynomials $q^{\prime }_{m,m+j} $ having degree $m+j-1$ by $P_{0,0},...,P_{0,m+j-1},P_{0,m+j},...,P_{0,M-1}:$
\begin{eqnarray}
q^{\prime }_{m,m+j} \left( \frac{s-a}{b-a} \right) = \sum _{l=0}^{M-1} \zeta_{j,l}^{m} P_{0,l} \left( \frac{s-a}{b-a} \right)
= \sum _{l=0}^{M-1} \zeta_{j,l}^{m} p_{0,l} (s), \label{gy100}
\end{eqnarray}
where $\zeta_{j,l}^{m} =0,$ if $m+j \leq l \leq M-1,$ and $\nu_m +m \leq M.$ Thus
\begin{eqnarray*}
\theta _{m,j} = q_{m,m+j} (1)f(b)-q_{m,m+j} (0)f(a) - \frac{1}{b-a} \sum _{l=0}^{M-1} \zeta _{j,l}^{m} \phi _l.
\end{eqnarray*}
By straightforward calculation
\begin{eqnarray*}
q_{m,m+j} (1)=1, \; q_{m,m+j} (0)= \left\{ 
                                           \begin{array}{ccc}
                                             (-1)^j , & \mbox{if} & m=0, \\
                                             0, & \mbox{if} & m > 0. \\
                                           \end{array}
                                         \right.
\end{eqnarray*}
Set now
\begin{eqnarray}
\underline{\zeta} _{j}^{o} &=& \left( 1, \; (-1)^{j+1} , \; -\zeta_{j,0}^{0},\ldots , -\zeta_{j,M-1}^{0} \right),\hspace{0.5cm}
 (m=0), \label{gy110} \\
\underline{\zeta} _{j}^{m} &=& \left( 1, \; 0, \; -\zeta_{j,0}^m,\ldots , -\zeta_{j,M-1}^m \right), \hspace{1.7cm}
 (m>0) \label{gy120}
\end{eqnarray}
 for $0 \leq j \leq \nu _m,$ and
\begin{eqnarray}
\mathcal{Z}_m &=&  \left[ (\underline{\zeta} _{0}^{m})^T,..., (\underline{\zeta} _{\nu_m}^{m})^T \right]^T \in
\mathbf{R}^{(\nu_m +1)\times (M+2)}. \label{gy50}
\end{eqnarray}
Estimation (\ref{e12}) can be obtained by direct substitution  taking into account (\ref{gy30}) and (\ref{e3b}). $\Box$

\begin{rem}\label{rem:3}
Paper \cite{zha15} gives the lower bound for functionals applied to derivatives of functions for several special cases together with
comparisons with previously published estimations, therefore we refer the reader for discussions to \cite{zha15}.
We only mention that the relation between Lemma \ref{lem:31} and Lemma 6 of \cite{leejfi15} is analogous to the one pointed out in Remark \ref{rem:2}.
Moreover neither \cite{leejfi15} nor \cite{zha15} derive any estimation for functionals applied to derivative of functions relying to polynomials of degree higher than $1.$
\end{rem}

\section{Stability analysis of continuous delayed systems}

Consider equation (\ref{x1}). Let  $M> 0, \; m_1 \geq 0, \; m_2 \geq 1$
be given integers.
Let $x_t (s)=x(t+s)$ be the solution of (\ref{x1}), and let
$\phi _j (t)  $   and $\Phi_M (t)$  be defined for function $f=x_t$ as before  with
$\phi_j (t) = \int _{-\tau}^{0} p_{0,j} (s) x_t (s) ds,$  and
$\Phi_M (t)= \mbox{col} \left\{  \phi_0 (t), ...,  \phi_{M-1} (t) \right\}$,
 Set furthermore
\begin{eqnarray*}
\widetilde{x} (t) = \mbox{col} \left\{ x(t), \Phi _M(t) \right\} , \;
\widetilde{\Phi}_M (t) = \mbox{\emph{col}} \left\{ x(t), x(t- \tau ), \frac{1}{\tau} \Phi _M(t) \right\}.
\end{eqnarray*}
Consider the LKF candidate
\begin{eqnarray}
V(x_t , \dot{x}_t ) = V_1 (x_t) + V_2 (x_t) + V_3 (\dot{x}_t),  \label{e15}
\end{eqnarray}

\vspace{-0.5cm}
where
\begin{eqnarray}
V_1 (x_t) &=&  \widetilde{x} (t)^T P \widetilde{x} (t), \hspace{5cm} P\in  \mathbf{S}_{n_x(M+1)},   \label{gy60} \\
V_2 (x_t) &=&  
 \sum _{j=0}^{m_1} \int _{-\tau}^{0}    \left(  \frac{s+ \tau}{\tau }  \right)^j     x_t (s)^T Q_j x_t (s)ds ,
 \hspace{0.8cm}Q_j \in \mathbf{S}_{n_x}^{+}, \hspace{2mm} j=0,...,m_1, \label{gy61} \\
V_3 (\dot{x}_t) &=&
\tau \sum _{j=1}^{m_2} \int _{-\tau}^{0} \left( \frac{s+ \tau}{\tau }\right)^j \dot{x}_t (s)^T R_j \dot{x}_t (s)ds ,
\hspace{0.5cm} R_j \in \mathbf{S}_{n_x}^{+}, \hspace{2mm} j=1,...,m_2.
\label{gy62}
\end{eqnarray}
We note that $V_2$ and $V_3$ can also be written as multiple integrals (c.f. (\ref{e400}), (\ref{e5})).

\begin{thm} \label{Th:13}
Let  $M> 0, \; m_1 \geq 0, \; m_2 \geq 1$ and $\nu_{1,j}\geq 0, \; (j=0, \ldots , m_1)$,  $\nu_{2,j} \geq 0, \; (j=0, \ldots , m_2)$
be given integers satisfying the inequalities $m_1+\nu_{1,j}<M,$ $m_2+\nu_{2,j} \leq M,$ for all $j$.
System (\ref{x1}) is asymptotically stable, if there are matrices $P \in \mathbf{S}_{n_x(M+1)}$, $Q_j \in \mathbf{S}_{n_x}^{+}$, $j=0,...,m_1$ and
$R_j \in \mathbf{S}_{n_x}^{+}$, $j=1,...,m_2$ such that the LMIs
\begin{eqnarray}
\Psi_{M,m_1}^0(\tau) >0, 
 \hspace{0.5cm}\Psi _{M}^{1}(\tau) + \Psi _{M,m_1}^{2} + \Psi _{M,m_2}^{3,1}(\tau) - \Psi _{M,m_2}
^{3,2}(\tau) < 0   \label{e20}
\end{eqnarray}
hold true, where

\begin{eqnarray}
\Psi_{M,m_1}^0 (\tau) &=& \tau P +  \sum _{j=0}^{m_1} \mbox{diag}
\left\{ 0, \left( \Xi_j \otimes I \right)^T \mathcal{Q}_{j}^{(j)}  \left( \Xi_j \otimes I \right)\right\}, \label{gy55} \\
 \Psi _{M}^{1}(\tau) &=& \Gamma_{M}^{T} P \Lambda_M +  \Lambda_{M}^{T} P \Gamma_M, \label{gy580} \\
\Psi _{M,m_1}^{2} &=& \mbox{diag} \left\{  \sum_{j=0}^{m_1} Q_j , \; -Q_0 , \; - \sum_{j=1}^{m_1} j
\left( \Xi _{j-1} \otimes I \right)^T \mathcal{Q}_{j-1}^{(j)} \left( \Xi _{j-1} \otimes I \right) \right\} ,\label{gy56} \\
\Psi _{M,m_2}^{3,1}(\tau) &=& \tau \mathcal{A}^T \sum _{j=1}^{m_2} R_j \mathcal{A} , \label{gy57}\\
\Psi _{M,m_2}^{3,2}(\tau) &=& \sum_{j=1}^{m_2} j
\left( \mathcal{Z}_{j-1} \otimes I \right)^T \mathcal{R}_{j-1}^{(j)} \left( \mathcal{Z}_{j-1} \otimes I \right),\label{gy58}
\end{eqnarray}
matrices $\Xi_k$ and $\mathcal{Z}_{k}$ are given by (\ref{gy20}) and (\ref{gy50}) with $\nu_{1,k}$ and $\nu_{2,k}$, respectively,
\begin{eqnarray}
\mathcal{Q}_{j}^{(k)} &=& \mbox{diag}  \left\{ (j+1)Q_k, \ (j+3)Q_k, \ ...\ , \ (j+(2M-1))Q_k \right\}, \label{gy81}\\
\mathcal{A}&=& \left( A, \; A_{d_1} , \; \tau A_{d_2} , 0 ,..., 0 \right)\in \mathbf{R}^{n_x \times n_x(M+2)}, 
 \label{gy82}\\
\Lambda _M &=&  \begin{bmatrix}
     		\mathcal{A} \\
      		\widetilde{\mathcal{L}}_0 \otimes I \\
    	\end{bmatrix}  , \;
\Gamma _M =  \begin{bmatrix}
     	1 & 0 & 0\\
        0 & 0 & \tau I_{M}\\
    	   \end{bmatrix} \otimes I, \;
     \label{gy83}\\
 \widetilde{\mathcal{L}}_0 &=&   \left[ (\underline{\zeta} _{0}^{0})^T,...,
(\underline{\zeta} _{M-1}^{0})^T \right]^T
  \label{gy840} \\
\mathcal{R}_{j-1}^{(j)}&=&\mbox{diag} \left\{ j  R_j, \ (j+2)R_j, \ ... \ , \ (j+2\nu_{j-1}) R_j \right\}, \label{gy84}
\end{eqnarray}
 where
$\underline{\zeta} _{j}^{0}$ is given by (\ref{gy110}),  and $0$ denotes zero matrices of compatible size.
\end{thm}

\textbf{Proof.}
We shall prove first the existence of a $\mu _1>0$ such that $V(x_t,\dot{x}_t) \geq \mu _1 \left\| x(t) \right\|$.
Consider the term of $V_2$ with $j=0$. Applying estimation
(\ref{e10}) with  $\nu_0=M-1$ and $\Xi_0=I$ we obtain
\begin{eqnarray}
\int _{-\tau}^{0} x_t (s)^T Q_0 x_t (s)ds &\geq&
 \tau ^{-1} \widetilde{x} (t)^T \mbox{diag} \left\{ 0 , \mathcal{Q}_{0}^{(0)} \right\}  \widetilde{x} (t).   \label{e17}
\end{eqnarray}
If $m_1 > 0$, apply (\ref{e10}) to the terms of $V_2$ with $j>0$. We obtain
\begin{eqnarray}
&&\int _{-\tau}^{0} \left(\frac{s+\tau}{\tau}\right)^j x_t (s)^T Q_j x_t (s)ds 
\geq -\tau ^{-1} \widetilde{x} (t)^T \mbox{diag} \left\{ 0_n , \left( \Xi _j \otimes I \right)^T \mathcal{Q}_{j}^{(j)}
\left( \Xi _j \otimes I \right) \right\}  \widetilde{x} (t).  \label{e18}
\end{eqnarray}
It follows from (\ref{e17}) and (\ref{e18}) that
\begin{eqnarray}
V_1 (x_t) + V_2 (t) \geq \widetilde{x} (t)^T \left( P+ \tau ^{-1} \sum _{j=0}^{m-1} \mbox{diag} \left\{ 0,
\left( \Xi _j \otimes I \right)^T \mathcal{Q}_{j}^{(j)} \left( \Xi _j \otimes I \right) \right\}  \right)  \widetilde{x} (t).   \label{e19}
\end{eqnarray}
Since $V_3 (\dot{x}_t) \geq 0,$ the existence of an appropriate $\mu _1$ follows from (\ref{e19}).

We shall prove next the negativity of $\frac{d}{dt} V(x_t , \dot{x}_t).$ Introduce the notation $\overline{V}_i (t) = V_i (x_t)$, $i=1,2,$ and
$\overline{V}_3 (t) = V_3 (\dot{x}_t)$, where $x_t$ is the solution of system (\ref{x1}).
The derivative of the first term of (\ref{e15}) is
\begin{eqnarray*}
\dot{\overline{V}}_1 (t)= \dot{\widetilde{x}} (t)^T P \widetilde{x}(t) + \widetilde{x}^T (t)P \dot{\widetilde{x}}(t),
\end{eqnarray*}
where $\dot{\widetilde{x}}(t)= \mbox{col} \left\{ \dot{x}(t) , \frac{d}{dt} \phi_0 (t),...,\frac{d}{dt} \phi_{M-1} (t)  \right\}$. The derivatives of $\phi_j$s
can be obtained by integration by parts:
\begin{equation}
\frac{d}{dt} \phi_j (t)= \int _{-\tau}^{0} p_{0,j} (s) \dot{x}_t (s) ds 
 = p_{0,j}(0)x(t)-p_{0,j}(-\tau)x(t-\tau)-  \frac{1}{\tau} \int _{-\tau}^{0} P_{0,j}^{\prime} \left( \frac{s+ \tau}{\tau} \right)
x_t (s)ds.     \label{bet1}
\end{equation}
In consistence with (\ref{gy01}), (\ref{e900}) and (\ref{gy20}), it follows from from (\ref{bet1}) that
\begin{eqnarray*}
\frac{d}{dt} \phi_j (t)=x(t)-(-1)^j x(t-\tau)- \sum _{l=0}^{M-1} \zeta _{j,l}^{0} \frac{1}{\tau} \phi _l(t) .
\end{eqnarray*}
Therefore, $\dot{\widetilde{x}} (t)= \Lambda_M \widetilde{\Phi}_M (t).$ On the other hand, $\widetilde{x} (t) = \Gamma _M \widetilde{\Phi}_M (t),$ thus we obtain
\begin{eqnarray}
\dot{\overline{V}}_1 (t) = \widetilde{\Phi}_M (t)^T \Psi _{M}^{1}(\tau) \widetilde{\Phi}_M (t).                        \label{e21}
\end{eqnarray}

The derivative of the first term of $V_2$ is
\begin{eqnarray}
\frac{d}{dt} \int _{-\tau}^{0} x_t (s)^T Q_0 x_t (s)ds = x(t)^T Q_0 x(t)-x(t-\tau)^T Q_0 x(t-\tau),   \label{e23}
\end{eqnarray}
while the derivatives of the terms of $V_2$ corresponding to $j \geq 1$ can be obtained as
\begin{align}
\frac{d}{dt} \int _{-\tau}^{0} \left(  \frac{s+\tau}{\tau} \right)^j  x_t (s)^T Q_j x_t (s)ds &=
 x(t)^T Q_j x(t)-\frac{j}{\tau} \int_{t-\tau}^{t} \left(\frac{s-t+\tau}{\tau} \right)^{j-1}  x(s)^T Q_j x(s)ds     \nonumber  \\
&= x(t)^T Q_j x(t)-\frac{j}{\tau} J_{Q_j , j-1,  -\tau , 0} ( x_t ).              \label{e24}
\end{align}
Employing Lemma \ref{lem:21}, we obtain from (\ref{e24}) that
\begin{eqnarray}
J_{Q_j ,j-1, -\tau , 0} ( x_t ) \geq \tau ^{-1} \Phi _M (t) ^T
\left( \Xi _{j-1} \otimes I \right)^T \mathcal{Q}_{j-1}^{(j)} \left( \Xi _{j-1} \otimes I \right) \Phi _M (t).   \label{e25}
\end{eqnarray}
It follows from (\ref{e25}) that
\begin{align*}
&\frac{d}{dt} \int _{-\tau}^{0} \left(  \frac{s+\tau}{\tau} \right)^j x_t (s)^T Q_j x_t (s)ds   \\ 
&\hspace{1.5cm} \geq  x(t)^T Q_j x(t) - j \frac{1}{\tau} \Phi _M (t)^T
\left( \Xi _{j-1} \otimes I \right)^T \mathcal{Q}_{j-1}^{(j)} \left( \Xi _{j-1} \otimes I \right) \frac{1}{\tau} \Phi_M (t),
\end{align*}
which means that
\begin{eqnarray}
\dot{\overline{V}}_2 (t) \leq \widetilde{\Phi}_M ^T \Psi _{M,m_1}^{2} \widetilde{\Phi}_M.                           \label{e26}
\end{eqnarray}
Now compute the derivative of $\overline{V}_3(t)$. We obtain
\begin{eqnarray}
\dot{\overline{V}}_3 (t)&=& \tau \sum _{j=1}^{m_2} \frac{d}{dt}  \int_{-\tau}^{0}\left(  \frac{s+\tau}{\tau} \right)^j
 \dot{x}_t (s)^T R_j \dot{x}_t (s)  ds= \nonumber  \\
&=& \tau \sum_{j=1}^{m_2} \left\{
\dot{x} (t)^T R_j \dot{x} (t) - \frac{j}{\tau} J_{j-1} (R_j, \dot{x}_t, -\tau, 0)
\right\} .                                   \label{e27x}
\end{eqnarray}
Applying now Lemma \ref{lem:31}, it follows  that
\begin{eqnarray}
J_{j-1} (R_j, \dot{x}_t, -\tau, 0)  \geq \tau^{-1} \widetilde{\Phi}_M (t)^T
\left( \mathcal{Z} _{j-1} \otimes I \right)^T \mathcal{R}_{j-1}^{(j)} \left( \mathcal{Z} _{j-1} \otimes I \right) \widetilde{\Phi}_M (t),
    \label{e27}
\end{eqnarray}
where $\mathcal{R}_{j-1}^{(j)}$ is given by (\ref{gy84}). 
From  (\ref{e27x}) and (\ref{e27}) we obtain
\begin{eqnarray}
\dot{\overline{V}}_3 (t) \leq \dot{x}(t)^T \left( \tau \sum_{j=1}^{m_2} R_j \right) \dot{x}(t)-
\tau^{-1} \widetilde{\Phi}_M (t)^T  \Psi _{M,m_2}^{3,2}(\tau)  \widetilde{\Phi}_M (t).     \label{e28}
\end{eqnarray}
Since $\dot{x} (t)= \mathcal{A} \widetilde{\Phi}_M (t),$  (\ref{e28}) implies that
\begin{eqnarray}
\dot{\overline{V}}_3 (t) \leq   \widetilde{\Phi}_M (t)^T  \left(\Psi _{M}^{3,1}(\tau)  
-
\Psi _{M,m_2}^{3,2}(\tau) \right)  \widetilde{\Phi}_M (t).                       \label{e29}
\end{eqnarray}
The statement of the theorem follows from  (\ref{e20}),  (\ref{e21}), (\ref{e26}) and (\ref{e29}) using the standard Lyapunov-Krasovskii Theorem (see e.g. \cite{frid2014}).    $\Box$

\begin{rem}
 If $A_{d_2}=0$ is considered in (\ref{e1}), Theorem \ref{Th:13} with $M=N, \; m_1=0, \; m_2=1 $ gives back Theorem 5 of \cite{seur14b} (apart from a multiplier $\tau$ (i.e. $h$) in the derivative of $V_3.$)
\end{rem}
\begin{rem} \emph{Delay range stability.}
An analogous stability result can be proven, if $A_d =0$ and $\tau$ is supposed to be an \textit{unknown constant}, but for which a lower and an upper bound is known, i.e.
$\underline{\tau} \leq \tau \leq \overline{\tau}$ for some given $\underline{\tau}$ and $\overline{\tau}.$ One can see that
$\Psi_{M,m_1}^{0}(\tau)$ is affine in $\tau ,$ and $\Psi_{M,m_1}^{0}(\tau)>0$ holds true for all
$\tau \in [\underline{\tau},\overline{\tau}], $ provided that $\Psi_{M,m_1}^{0}(\overline{\tau})>0.$ Moreover
$\Psi_{M,m_1}^{2}$ does not depend on $\tau ,$  while in this case, $\Psi_{M}^{1}(\tau)$ is affine in $\tau ,$ as well.
One can modify the definition (\ref{gy62}) of $V_3(\dot{x}_t)$ by taking the multiplier $\tau ^2$ in front of the summation instead of $\tau,$
then $\tau ^{-1}$ disappears from $\Psi_{M,m_2}^{3,2}(\tau).$
Apply  Schur complements to (\ref{e20}) with respect to the new $\Psi_{M}^{3,1}(\tau)$ and a congruence transformation, then we obtain
\begin{eqnarray}
\overline{\Psi}(\tau) =
   \begin{bmatrix}
       \Psi _{M}^{1}(\tau) + \Psi _{M,m_1}^{2} - \Psi _{M,m_2}^{3,2}    &   \tau \mathcal{A}^{T}  \sum_{j=1}^{m} R_j \\
        *  &    -\sum_{j=1}^{m} R_j \\
      \end{bmatrix} < 0 .          \label{e31}
\end{eqnarray}
The matrix valued function $\overline{\Psi}(\tau)$ is affine in $\tau$,
which means that it is enough to require the fulfillment of the inequality (\ref{e31}) at the endpoints, i.e.
the LMIs  $\Psi _{M,m_1}^{0}(\overline{\tau})>0,$ $\overline{\Psi}(\underline{\tau}) <0$ and
$\overline{\Psi}(\overline{\tau}) <0$  have to hold true.
\end{rem}

\section{Hierarchy of the LMI stability conditions}

This section is devoted to the comparison of the stability conditions obtained in the previous section for different parameters.
We observe that parameter $M$ determines the size of matries $P$  and $\widetilde{\mathcal{L}}_0$, the number of columns of $\Xi_j$ and $\mathcal{Z}_k$,
while the number of rows of $\Xi_j$ and $\mathcal{Z}_k$ is $\nu_{1,j}$ and $\nu_{2,k}$. The number of matrices $Q_j$s and $R_k$s is $m_1$ and $m_2$. \emph{The aim is to show that the LMI conditions can be arranged into a hierarchy table provided that the parameters are chosen to satisfy the following condition.}
\begin{eqnarray}
\begin{array}{lll}
   M\geq 1, & m_1 =m, & m_2 = m+1, \\
  \nu_{1,j}=
  M-j-1, & \nu_{2,j}=\nu_{1,j}+1, & j=0,1,...,m .
\end{array}
  \label{H1}
\end{eqnarray}
We shall refer to the LMI condition (\ref{e20}) with parameters satisfying (\ref{H1}) as $\mathcal{L}_{M,m}(\tau)$.
\begin{defi}
Let the pairs $(M,m)$ and $(\hat{M},\hat{m})$ be given. We will say that
\emph{$\mathcal{L}_{\hat{M},\hat{m}} $ outperforms $\mathcal{L}_{M,m} $,} if, for every $\tau$ for which
$\mathcal{L}_{M,m} (\tau)$ has a feasible solution, $\mathcal{L}_{\hat{M},\hat{m}} (\tau)$ has a feasible solution, too.
This is denoted by $\mathcal{L}_{M,m} \prec \mathcal{L}_{\hat{M},\hat{m}}.$
\end{defi}
We will show that the parametric family of $\mathcal{L}_{M,m}$ is ordered according to both parameters.


\begin{thm} \label{Th:4}
Let the integer parameters satisfy (\ref{H1}).  Then
\begin{eqnarray}
  \mathcal{L}_{M,m} &\prec& \mathcal{L}_{M+1,m},  \label{H10}\\
  \mathcal{L}_{M,m} &\prec& \mathcal{L}_{M,m+1}.  \label{H11}
\end{eqnarray}
\end{thm}
\textbf{Proof.} \emph{Part 1.}
First we show (\ref{H10}).

Let matrices $P$, $Q_{0},...,Q_{m-1}$ and
$R_{1},...,R_{m}$ denote a feasible solution of $\mathcal{L}_{M,m} (\tau)$ for some fixed $\tau .$ We seek the solution of $\mathcal{L}_{M+1,m} (\tau)$
in the form of
\begin{equation} \label{H20}
\hat{P}= \begin{bmatrix}
               P  &   0   \\
                 0    &   \varepsilon I   \\
           \end{bmatrix} ,
\hspace{2mm} \hat{Q}_{i} = Q_{i} , \; (i=0,...,m-1), \hspace{2mm} \hat{ R}_{j} = R_{j}, \; (j=1,...,m )
\end{equation}
for some positive constant $\varepsilon .$ In what follows, we shall denote matrices that belong to $(M+1,m)$ analogously by putting "hat" over them.

We show first that inequality $\Psi_{M,m}^{0} (\tau)>0 $ implies  $ \Psi_{M+1,m}^{0} (\tau)$
independently of
$\varepsilon>0$. In fact,
 matrix $\hat{\Xi}_{j}$ is obtained by adding a new row and a new column to
$\Xi_{j}$, i.e
\begin{equation*}
\hat{\Xi}_{j} = \begin{bmatrix}
                    \Xi_{j}  &  0    \\
                     \underline{\xi}_{M-j}^{j,1} &  \xi_{M-j,M}^{j}      \\
                  \end{bmatrix} ,
\end{equation*}
where $\underline{\xi}_{M-j}^{j}$ is partitioned as
$\underline{\xi}_{M-j}^{j} = \left( \underline{\xi}_{M-j}^{j,1}, \xi_{M-j,M}^{j} \right).$
Thanks to the structure of the matrices,
we obtain by standard algebra
\begin{eqnarray}
&& \hspace{-0.8cm}
\Psi_{M+1,m}^{0}(\tau)=
                \begin{bmatrix}
                    \Psi_{M,m}^{0} (\tau)&  0    \\
                     0             &  \varepsilon I  \\
                  \end{bmatrix}
+ \mbox{PSDTs,}
               \label{e32}
\end{eqnarray}
where PSDTs stands for positive semidefinite terms, therefore the statement follows.

 Next we express $\Psi _{M+1}^{1}(\tau)$ by $\Psi _{M}^{1} (\tau).$
By the definition of $\Psi _{M+1}^{1}(\tau)$, we have
\begin{eqnarray*}
\Psi_{M+1}^{1}(\tau) &=& \Gamma_{M+1}^{T} \hat{P} \Lambda_{M+1} +  \Lambda_{M+1}^{T} \hat{P} \Gamma_{M+1}, \hspace{1cm} \mbox{with}   \\
\Gamma _{M+1} &=& \mbox{diag} \left\{  \Gamma_M , \tau I \right\} , \;
\Lambda_{M+1} =
                  \begin{bmatrix}
                    \Lambda_{M} &  0    \\
                    \underline{\zeta}_{M}^{0,1} \otimes I &  \underline{\zeta}_{M}^{0,2} I  \\
                  \end{bmatrix} , \\
\end{eqnarray*}
where $\underline{\zeta}_{M}^{0} $ is again partitioned as $\underline{\zeta}_{M}^{0} = \left(\underline{\zeta}_{M}^{0,1}, \underline{\zeta}_{M}^{0,2}\right)$ with $\underline{\zeta}_{M}^{0,2}= -\zeta_{M,M}^{0}.$
Using (\ref{H20}), by standard algebra we obtain
\begin{equation}
\Psi_{M+1}^{1}(\tau) =
               \begin{bmatrix}
                    \Psi_{M}^{1}(\tau) &  0    \\
                    0 &  0  \\
                  \end{bmatrix}
+ \varepsilon \tau \left(
               \begin{bmatrix}
                    0 &  0    \\
                    \underline{\zeta}_{M}^{0,1} \otimes I &  \underline{\zeta}_{M}^{0,2} I  \\
                  \end{bmatrix}  +
                \begin{bmatrix}
                    0 &  \left( \underline{\zeta}_{M}^{0,1} \otimes I \right)^T    \\
                     0 &  \underline{\zeta}_{M}^{0,2} I \\
                  \end{bmatrix}
\right) .                          \label{e34}
\end{equation}

Express now $\Psi _{M+1,m}^{2}$ by $\Psi _{M,m}^{2} .$ Since
\begin{equation*}
\hat{\Xi}_{j-1}=
               \begin{bmatrix}
                    {\Xi}_{j-1} &  0    \\
                    \underline{\xi}_{M+1-j}^{j-1,1} &  \xi _{M+1-j,M}^{j-1}  \\
                  \end{bmatrix} , \;
\hat{\mathcal{Q}}_{j-1}^{(j)} = \mbox{diag} \left\{  {\mathcal{Q}}_{j-1}^{(j)} ,c_2(M,j)Q_j \right\}   ,
\end{equation*}
where $c_2(M,j)=2M-j+2$, we obtain
\begin{eqnarray*}
\left(  \hat{\Xi}_{j-1} \otimes I  \right)^T \hat{\mathcal{Q}}_{j-1}^{(j)} \left( \hat{ \Xi} _{j-1} \otimes I  \right)& =&
               \begin{bmatrix}
                    \left(  \Xi_{j-1} \otimes I  \right)^T \mathcal{Q}_{j-1}^{(j)}
                    \left(  \Xi_{j-1} \otimes I  \right) &  0    \\
                    0 &  0  \\
                  \end{bmatrix} + \mbox{PSDTs,}  
\end{eqnarray*}
therefore
\begin{eqnarray}
&&\Psi _{M+1,m}^2(\tau) \leq
                \begin{bmatrix}
                    \Psi _{M,m}^2(\tau) &  0    \\
                    0 &   0   \\
                 \end{bmatrix} .
 \label{e35}
\end{eqnarray}
Express now $\Psi _{M+1,m}^{3,1}(\tau)$ by $\Psi _{M,m}^{3,1}(\tau) .$
Since $\hat{\mathcal{A}}= \left(  \mathcal{A},0 \right),$ we obtain
\begin{equation}
\Psi _{M+1,m}^{3,1}(\tau) =
                \begin{bmatrix}
                    \Psi _{M,m}^{3,1}(\tau) &  0    \\
                    0 &   0   \\
                  \end{bmatrix}  .            \label{e36}
\end{equation}
 Express $\Psi _{M+1,m}^{3,2}(\tau)$ by $\Psi _{M,m}^{3,2}(\tau) .$
Since
\begin{equation*}
\hat{\mathcal{Z}}_{j-1}^{(j)}=
                \begin{bmatrix}
                    \mathcal{Z}_{j-1}^{(j)} &  0    \\
                    \underline{\zeta}_{M+1-j}^{j-1,1} &  \underline{\zeta}_{M+1-j}^{j-1,2}  
                  \end{bmatrix} , \;
\hat{\mathcal{R}}_{j-1}^{(j)} = \mbox{diag} \left\{  \mathcal{R}_{j-1}^{(j)} ,c_2(M,j) R_j \right\}   ,
\end{equation*}
where $\underline{\zeta}_{M+1-j}^{j-1}= \left[ \underline{\zeta}_{M+1-j}^{j-1,1},\underline{\zeta}_{M+1-j}^{j-1,2}\right]$ and
$\underline{\zeta}_{M+1-j}^{j-1,2}=-{\zeta}_{M+1-j,M}^{j-1}$, we obtain
\begin{eqnarray}
\Psi _{M+1,m}^{3,2}(\tau) &=& \frac{1}{\tau}\sum_{j=1}^{m}  j \left(  \hat{\mathcal{Z}}_{j-1}^{(j)} \otimes I  \right)^T \hat{\mathcal{R}}_{j-1}^{(j)} \left(  \hat{\mathcal{Z}}_{j-1}^{(j)} \otimes I  \right)=   \nonumber    \\
&\geq&\frac{1}{\tau} \sum_{j=1}^{m}  j
  \begin{bmatrix}
                    \left(  \mathcal{Z}_{j-1}^{(j)} \otimes I  \right)^T \mathcal{R}_{j-1}^{(j)}
                    \left(  \mathcal{Z}_{j-1}^{(j)} \otimes I  \right) &  0    \\
                    0 &   0   \\
                  \end{bmatrix} \nonumber\\
&&
+\frac{c_2(M,1)}{\tau}
                 \begin{bmatrix}
                    0 &  \left(  \underline{\zeta}_{M}^{0,1} \otimes I  \right)^T     \\
                    0 &   -\zeta_{M,M}^{0}  I    \\
                   \end{bmatrix}
                   \begin{bmatrix}
                   0 & 0 \\
                   0&R_1 \\
                   \end{bmatrix}
                   \begin{bmatrix}
                    0 &  \left(  \underline{\zeta}_{M}^{0,1} \otimes I  \right)^T     \\
                    0 &   -\zeta_{M,M}^{0}  I    \\
                   \end{bmatrix}^T 
,                                                          \label{e37}
\end{eqnarray}
where we omitted several positive semidefinite terms on the right hand side of (\ref{e37}).

Finally we show that $\overline{\Psi}_{M+1}(\tau)= \Psi_{M+1}^{1}(\tau)+ \Psi_{M+1,m}^{2}+ \Psi_{M+1}^{3,1}(\tau)
 - \Psi_{M+1}^{3,2} (\tau)< 0.$
 Applying (\ref{e34})-(\ref{e37}) we obtain
\begin{eqnarray}
\overline{\Psi}_{M+1}(\tau) &\leq&
  \begin{bmatrix}
        I & \left(  \underline{\zeta}_{M}^{0,1} \otimes I  \right)^T  \\
        0 & -\zeta_{M,M}^{0}  I \\
      \end{bmatrix}
    \begin{bmatrix}
        \overline{\Psi}_{M}(\tau) & 0 \\
        0 & -\frac{2M+1}{\tau}R_1 \\
      \end{bmatrix}
      \begin{bmatrix}
                    I &  0       \\
                    \left(  \underline{\zeta}_{M}^{0,1} \otimes I  \right) &   -\zeta_{M,M}^{0}  I   \\
                  \end{bmatrix}
    \nonumber   \\
    & +& \varepsilon \tau
    \begin{bmatrix}
        0 & \left(  \underline{\zeta}_{M}^{0,1} \otimes I  \right)^T  \\
        \underline{\zeta}_{M}^{0,1} \otimes I & 2 \underline{\zeta}_{M}^{0,2}I \\
      \end{bmatrix}
 \label{e38}
\end{eqnarray}
Since $\mbox{diag} \left\{ \overline{\Psi}_{M}(\tau), -\frac{2M+1}{\tau}R_1  \right\} <0,$ and matrix
$\begin{bmatrix}
                    I &  0       \\
                    \left(  \underline{\zeta}_{M}^{0,1} \otimes I  \right) &   -\zeta_{M,M}^{0}  I   \\
                  \end{bmatrix}$
                  is non-singular, there exists a positive constant $\nu_1$ such that
\begin{equation*}
\begin{bmatrix}
        I & \left(  \underline{\zeta}_{M}^{0,1} \otimes I  \right)^T  \\
        0 & -\zeta_{M,M}^{0}  I \\
      \end{bmatrix}
    \begin{bmatrix}
        \overline{\Psi}_{M}(\tau) & 0 \\
        0 & -\frac{2M+1}{\tau}R_1 \\
      \end{bmatrix}
    \begin{bmatrix}
    I &  0       \\
                    \left(  \underline{\zeta}_{M}^{0,1} \otimes I  \right) &   -\zeta_{M,M}^{0}  I   \\
                  \end{bmatrix}
     < -\nu_1 I,
\end{equation*}
and constant $\nu_2$ such that
\begin{equation*}
\tau
    \begin{bmatrix}
        0 & \left(  \underline{\xi}_{M}^{0,1} \otimes I  \right)^T  \\
        \underline{\xi}_{M}^{0,1} \otimes I & 2 \underline{\xi}_{M}^{0,2}I \\
      \end{bmatrix}  < \nu_2 I.
\end{equation*}
If $\varepsilon \nu_2 < \nu_1,$ then $\overline{\Psi}_{M+1}(\tau)<0$ holds true.

\emph{Part 2}. Secondly we show that $\mathcal{L}_{M,m} \prec \mathcal{L}_{M,m+1}.$
First we show the positivity of $\Psi_{M,m+1}^{0}(\tau).$ Suppose that $\Psi_{M,m}^{0}(\tau) >0$ with
the choice of $P$, $Q_0,...,Q_{m-1}.$ We seek matrix $Q_m $ as $Q_m = \varepsilon I$ with $\varepsilon>0.$ Then we obtain
\begin{eqnarray}
\Psi_{M,m+1}^{0}(\tau)
&=&  \Psi_{M,m}^{0}(\tau) +
\frac{\varepsilon}{\tau}
\mbox{diag} \left\{ 0,\left(  \Xi_{m} \otimes I  \right)^T  \mathcal{D}_{M,m} \left(  \Xi_{m} \otimes I  \right)   \right\},          \label{e40} \\
\mathcal{D}_{M,m}&=&\mbox{diag}\left\{ (m+1)I,(m+3)I,...,\left(2(M-m-1)+m+1\right)I \right\}.
\end{eqnarray}
The first term of the right hand side of (\ref{e40}) is positive, the second term is non-negative, therefore
$\Psi_{M,m+1}^{0}(\tau)>0 $ for any  $\varepsilon>0.$

Next we show that $\overline{\Psi}_{M,m+1}(\tau)<0$ has a feasible solution, provided that $\overline{\Psi}_{M,m}(\tau)<0$ has.
Suppose that $\overline{\Psi}_{M,m} (\tau)<-\nu _3 I$ with $P$, $Q_0,...,Q_{m-1}$ and $R_0,...,R_{m-1}$. We seek matrices $Q_m = \varepsilon I$
and $R_{m+1}= \varepsilon I$, where $\varepsilon >0$ has to be chosen. Matrix $\Psi_M^1(\tau)$ is unchanged.

Denote $E_1=\begin{bmatrix}
        I & 0 & \ldots 0  \\
      \end{bmatrix} \in \mathbf{R}^{n_x\times n_x(M+2)},$ then
\begin{equation}
\Psi _{M,m+1}^2 (\tau)= \Psi _{M,m}^2(\tau) + E_1^T Q_m E_1 -
(m+1)\left(  \Xi_{m} \otimes I  \right)^T  \mathcal{Q}_m \left(  \Xi_{m} \otimes I  \right)
   \label{e41}
\end{equation}
holds true.
On the one hand,
\begin{equation}
\Psi _{M,m+1}^{3,1}(\tau) = \tau \mathcal{A}^T \sum _{j=1}^{m} R_j \mathcal{A}+\tau \mathcal{A}^T  R_{m+1} \mathcal{A}=
\Psi _{M,m}^{3,1}(\tau)+\varepsilon \tau \mathcal{A}^T  \mathcal{A},   \label{e42}
\end{equation}
while
\begin{equation}
\Psi _{M,m+1}^{3,2}(\tau) = \Psi _{M,m}^{3,2}(\tau) + \frac{m+1}{\tau}\varepsilon
\left(  \mathcal{Z}_{m+1} \otimes I  \right)^T  \mathcal{D}_{M,m+1} \left(  \mathcal{Z}_{m+1} \otimes I  \right).   \label{e43}
\end{equation}
It follows from (\ref{e42})-(\ref{e43}) that
\begin{equation*}
\overline{\Psi}_{M,m+1} (\tau) = \overline{\Psi}_{M,m} (\tau) + \varepsilon \Omega_{M,m}(\tau)
\end{equation*}

\vspace{-0.5cm}
with
\begin{eqnarray*}
\Omega_{M,m}(\tau)&=&   E_1'  E_1 -
(m+1)\left(  \Xi_{m} \otimes I  \right)^T  \mathcal{D}_{M,m} \left(  \Xi_{m} \otimes I  \right)   \\
&&+\tau \mathcal{A}^T  \mathcal{A}-\frac{m+1}{\tau}
\left(  \mathcal{Z}_{m} \otimes I  \right)^T  \mathcal{D}_{M,m+1} \left(  \mathcal{Z}_{m} \otimes I  \right)
\end{eqnarray*}
Then there exists a constant $\nu_4$ such that $\Omega_{M,m}(\tau) \leq \nu_4 I.$ If $\varepsilon >0$ is small enough to satisfy
inequality $\varepsilon \nu_4 <\nu_3 ,$ then $\Psi_{M,m+1} (\tau)$ is negative.    $\Box$

\section{Numerical examples}

 In this section, we apply the proposed method to three benchmark examples that have been extensively used in the literature to compare the results.
The computations have been performed by using YALMIP \cite{yalmip} together with MATLAB.

\subsection{Some remarks on the implementation}

Assume that the integer parameters are chosed according to (\ref{H1}).
In order to implement LMIs (\ref{e20}) with (\ref{gy55})-(\ref{gy84}),  matrices $\widetilde{\mathcal{L}}_0, $ $\Xi _i$  and $\mathcal{Z}_j$ are to be produced.
This matrices can be computed employing the generalized Rodrigues formula (\ref{e3}) as follows. Let $X_K =\left(1, x, \ldots, x^K \right)^T $ and $\Pi_{m,K}(x)=\left(P_{m,0}(x), P_{m,1}(x), \ldots, P_{m,K}(x) \right)^T .$ Then $\Pi_{m,K}(x)=G(m,K)X_K$
 and $X_K= G(0,M-1)^{-1}\Pi_{0,K}(x),$  where $G(m,K) \in \mathbf{R}^{(K+1)\times (K+1)}$ with elements $G(m,K)_{1,1}=1,$
\begin{equation*}
  G(m,K)_{l+1,k+1}= (-1)^{l+k}\prod _{j=0}^{k-1}\frac{l-j}{k-j}\prod _{i=1}^{l}\frac{m+k+i}{i}, \;
  \mbox{ if } \hspace{3mm} \begin{array}{l}
                l=1,\ldots K-1, \\
                K=0,\ldots,l,
              \end{array}
  \end{equation*}
and $G(m,K)_{l,k}=0,$ if $k>l.$ By taking into account (\ref{e9}), (\ref{gy20}) and (\ref{H1}) we can see that
\begin{eqnarray*}
  \Xi_m &=& G(m,-\nu _m)
  \begin{bmatrix}
        0_{\nu _m+1,m} & I_{\nu _m+1} \\
      \end{bmatrix} G(0,M-1)^{-1}.
\end{eqnarray*}
Further, by taking into account (\ref{gy100}), (\ref{gy110}),
(\ref{gy50}) and (\ref{H1}) we can see that
\begin{eqnarray*}
\widetilde{\mathcal{L}}_0 &=&
\left[
  \begin{array}{ccc}
   \underline{\ell}_{M-1}^{(1)} & \underline{\ell}_{M-1}^{(2)} &  -L_0
  \end{array}
  \right]
   \\
   && \hspace{0.1cm} L_0 =\begin{bmatrix}   I_{M} \; 0\end{bmatrix} G(0,M) diag\left\{0,1,\ldots,M\right\}
    \begin{bmatrix} 0 \\ I_{M-1}\end{bmatrix} G(0,M-1)^{-1}, \\
\mathcal{Z}_0 &=&
  \left[
  \begin{array}{ccc}
   \underline{\ell}_{\nu_0}^{(1)} & \underline{\ell}_{\nu_0}^{(2)} &  -Z_0
  \end{array}
  \right]
   \\
   && \hspace{0.1cm} Z_0 = G(0,M-1) diag\left\{0,1,\ldots,M\right\}
    \begin{bmatrix} 0 \\ I_{M-1}\end{bmatrix} G(0,M-1)^{-1}, 
\end{eqnarray*}
\begin{eqnarray*}
\mathcal{Z}_m &=& \left[
  \begin{array}{ccc}
   \underline{\ell}_{\nu_m}^{(1)} & \underline{0}_{\nu_m} &-Z_m
    \end{array}
  \right]   \\
   && \hspace{0.1cm} Z_m=
   G(m,\nu _m+1)
   \left[
  \begin{array}{cc}
      0_{\nu_m+1, m-1} & D_{m,\nu _m}
     \end{array}
  \right]
     G(0,m+\nu_m)^{-1},
\end{eqnarray*}
where the vectors $\underline{\ell}_{k}^{(1)}, \underline{\ell}_{^k}^{(2)}, \underline{0}_{k} \in \mathbf{R}^{k+1}, $ are defined by
$\underline{\ell}_{k}^{(1)}=(1, \, \ldots, \, 1)^T,$ $\underline{0}_{k}=(0,\ldots,0)^T,$ $\underline{\ell}_{^k}^{(2)}=(-1, 1, \ldots, \pm 1)^T$ and
$D_{m,\nu _m}=\mbox{diag}\left\{m,  \ldots, m+\nu _m+1 \right\}.$
\subsection{Numerical experiments}
 Consider system (\ref{x1}) with coefficient matrices listed in Table 1. 
{\footnotesize
\begin{table}[!ht] \label{Tab:1}
\caption{{\footnotesize Systems used as illustrative examples}}
{\footnotesize
\begin{center}
\begin{tabular}{lcccc}
\toprule
Example  & $A$ & $A_{d_1}$ & $A_{d_2}$ & analytical bounds
\\
\hline

\vspace{-0.35cm}
 & & & & \\
\vspace{1mm}
1 & $\begin{bmatrix}-2 & 0 \\ 0 & -0.9  \end{bmatrix}$
   & $\begin{bmatrix}-1 & 0 \\-1 & -1  \end{bmatrix}$
   & $\begin{bmatrix} 0 & 0 \\0 & 0  \end{bmatrix}$
   & $\begin{array}{c}
        \underline{\tau}=0 \\
        \overline{\tau}\sim 6.17258
      \end{array}$  \\

   \hline

\vspace{-0.35cm}
 & & & & \\
\vspace{1mm}
2  & $\begin{bmatrix}0.2 & 0 \\ 0.2 & 0.1  \end{bmatrix}$
   & $\begin{bmatrix}0 & 0 \\0 & 0  \end{bmatrix}$
   & $\begin{bmatrix} -1 & 0 \\-1 & -1  \end{bmatrix}$
   & $\begin{array}{c}
        \underline{\tau}\sim0.2 \\
        \overline{\tau}\sim 2.04
      \end{array}$
     \\

   \hline

\vspace{-0.35cm}
 & & & & \\
\hspace{0.5mm}
3  & $\begin{bmatrix}0 & 1 \\ -2 & 0.1  \end{bmatrix}$
   & $\begin{bmatrix}0 & 0 \\1 & 0  \end{bmatrix}$
   & $\begin{bmatrix} 0 & 0 \\0 & 0  \end{bmatrix}$
   & $\begin{array}{c}
        \underline{\tau}\sim 0.1002 \\
        \overline{\tau}\sim 1.7178
      \end{array}$
    \\
  \bottomrule
 \end{tabular}  \\
 \end{center}
  \vskip-2mm
     }
   \end{table}
}

\emph{Example 1} is considered in numerous papers, among others, in \cite{seur14b,kim16,seur14}, where extensive comparisons with previously reported results are given. The results obtained by Theorem \ref{Th:13} are given in Table 2. 
The results obtained with $m=1,$ $M=1;2;3$ coincide with that of \cite{seur14b}, which are the best previously reported results we are aware of. We note that the same bounds of $\tau$ are obtained for $m=2$ and $m=3$ with the same values of $M,$ however an improvement is resulted in in the case of $m=4, \ M=3$ compared to $m=3, \ M=4.$ NoDV denotes in all tables the number of decision variables.

{\footnotesize
\begin{table}[!ht] \label{Tab:2}
\caption{{\footnotesize Delay bounds for Example 1 obtained by Theorem \ref{Th:13}}}
{\footnotesize
\begin{center}
\begin{tabular}{lccccc}
\toprule
$m$                &$1$ & $1$ & $1$ & $1$ & 4  \\
   \hline
$M$                & $1$     & $2$      & $3$          & $4$          &$3$ \\
   \hline
$\overline{\tau}$  &$6.05932$&$6.16893$ &$6.17250$     & $6.17258$ & $6.17258$\\
\hline
 NoDV \hspace{0.5cm}    & $16$ & $32$ & $48$     & $61$ &  $72$  \\ %
   \bottomrule
 \end{tabular}  \\
 \end{center}
  \vskip-2mm
    }
      \end{table}
      }
\emph{Example 2} is considered in many papers, among others, in \cite{park15,S-GAut13,sga15}, where extensive comparisons with previously reported results are given. In \cite{S-GAut13}, a delay bounding interval $[0.200, \ 1.877]$ was obtained with $16$ decision variables, while the authors of \cite{park15}
derived the delay bounding interval 
$[0.2000, \ 1.9504]$ from Theorem 1 with $59$ decision variables. In \cite{sga15}, the lower bound of the delay was found to be $0.2001$, while by Theorem \ref{Th:13}, we obtained the lower bound $0.20001.$ The upper bounds reported in \cite{sga15} and obtained by Theorem \ref{Th:13} are given in Table 3. The values of $N$ in \cite{sga15} and  $M$ in the present paper are related as $N=M+1.$   We note that the same upper bounds were obtained for different values of $m$ for a given value of $M.$ 
%
%
%
{\footnotesize
\begin{table}[!ht] \label{Tab:3}
\caption{{\footnotesize Delay upper bounds for Example 2 }} 
{\footnotesize
\begin{center}
\begin{tabular}{llccccc}
\toprule
Method &  $ M$\hspace{5mm}&    $ 1$ & $ 2$ & $3$ & $ 4$  & $5$ \\
  \hline
\cite{sga15}
&   $\overline{\tau}$  & $-$ & $1.58 \hspace{3.3mm}$ & $1.83\hspace{3.3mm}$ & $1.95\hspace{3.3mm}$     & $2.02\hspace{3.3mm}$ \\
  \hline
Theorem 1, $m=1$ \hspace{0.5cm}
&   $\overline{\tau}$  &$1.9419 $ & $2.0395$ &$2.0412$     & $2.0412$   & $2.0412$\\
   \hline
& NoDV \hspace{0.5cm}   & $16$ & $32$ & $48$     & $61$ &  $84$  \\ 
   \bottomrule
 \end{tabular}  \\
 \end{center}
  \vskip-2mm

    }
      \end{table}
      }

\emph{Example 3} is is also widely used for comparing the effectiveness of different methods. Here we shall mention a  recently published work in \cite{park15}, where comparisons with previously reported results are given, as well (see also \cite{trinh15}). The results obtained by \cite{park15} and by Theorem \ref{Th:13} are given in Table 4.  We note that the same upper bounds were obtained for different values of $m$ for a given value of $M$ for this example, too. 

\begin{table}[!ht] \label{Tab:4}
\caption{{\footnotesize Delay bounds for Example 3}}
{\footnotesize
\begin{center}
\begin{tabular}{lccccc}
\toprule
Method  \hspace{2cm}                &$ $  & $M$  & $\underline{\tau}$  & $\overline{\tau}$  &   NoDV\\ 
   \hline
 \cite{park15} Theorem 1                    & $ $ & $ $  & $   0.1002\hspace{2mm}$              & $1.5954\hspace{2mm}$                 & $59$\\
   \hline
 Theorem 1, $(m=1)$                         & $ $ & $1$  & $   0.10055$              & $1.5405\hspace{2mm}$                 &$16$\\
   \hline
                                  & $ $ & $2$  &   $0.10018 $            & $1.7122\hspace{2mm}$                 &$32$\\
   \hline
                                 & $ $ & $3$  &   $0.10017$            & $1.71799$                 & $48$\\
   \bottomrule
 \end{tabular}  \\
 \end{center}
  \vskip-2mm
     }
      \end{table}
\subsection{Discussion}
It can be seen that, in these examples, Theorem 1 yields better delay bounds than previously published methods except of Theorem 5 of \cite{seur14b} which is equivalent to Theorem 1, if it is applied with $m=1.$ It is worth noting that the better results are obtained with much smaller number of decision variables.

In these and in several other examples from the literature, on which we tested our approach, we observed that the improvement of the delay estimation is primarily due to the increase of the dimension of the extended state variable together with the improved lower bounds of Lemma \ref{lem:21} and Lemma \ref{lem:31}. This does not contradict to the reported improvements in the case of the application of triple, etc. integral terms in the LKF, since the applied lower estimations of the integrals of quadratic terms lead to the introduction of some extended state variables with increased dimension, as well. We emphasize that this remark is limited to the investigated examples, and the observed behavior may have the reason that the analytical bounds were rapidly reached up to 4-6 digits. Moreover, we note similarly to \cite{seur14b} that the formulated result does not establish any convergence to the analytical bounds.
\section{Conclusion}
In this paper, new multiple integral inequalities are derived based on certain hypergeometric-type orthogonal polynomials. These inequalities   are similar to that of \cite{zha15}, and they comprise the famous Jensen's and Wirtinger's inequalities, as well as the recently presented
Bessel-Legendre inequalities of \cite{seur14b} and the Wirtinger-based multiple-integral inequalities of \cite{leejfi15, park15}. Applying the obtained inequalities, a set of sufficient LMI stability conditions for linear continuous-time delay systems are derived. It was proven that these LMI conditions could be arranged into a bidirectional hierarchy establishing a rigorous theoretical basis for comparison of conservatism of the investigated methods.
Numerical examples confirm that the proposed method enhances the tolerable delay bounds.


\end{document}